\documentclass [12pt] {article}

\usepackage{latexsym}

\font\bb=msbm10 at 12pt

\newcommand{\D}{\mbox{{\bb D}}}
\newcommand{\C}{\mbox{{\bb C}}}

\everymath{\displaystyle}

\newenvironment{multi}
{ \begin{equation} \begin{array}{c} }
{ \end{array} \end{equation} }

\makeatletter

\@addtoreset{equation}{section}
\makeatother

\begin{document}

\title{New characterizations of ruled real hypersurfaces in complex projective space}
\author{Juan de Dios P\'{e}rez and David P\'{e}rez-L\'{o}pez}
\date{}
\maketitle

\begin{abstract} We consider real hypersurfaces $M$ in complex projective space equipped with both the Levi-Civita and generalized Tanaka-Webster connections. For any nonnull constant $k$ and any symmetric tensor field of type (1,1)  $L$  on $M$ we can define two tensor fields of type (1,2) on $M$, $L_F^{(k)}$ and $L_T^{(k)}$, related to both connections. We study the behaviour of the structure operator $\phi$ with respect to such tensor fields in the particular case of $L=A$, the shape operator of $M$, and obtain some new characterizations of ruled real hypersurfaces in complex projective space.
\end{abstract}



\section{Introduction.}
Let $\C P^m$, $m\geq 2$, be the complex projective space endowed with the Kaehlerian structure $(J,g)$, where $g$ is the Fubini-Study metric of constant holomorphic sectional curvature $4$. Let $M$ be a connected real hypersurface of $\C P^m$ without boundary, $g$ the restriction of the metric on $\C P^m$ to $M$ and $\nabla$ the Levi-Civita connection on $M$. Take a locally defined unit normal vector field $N$ on $M$ and denote by $\xi =-JN$. This is a tangent vector field to $M$ called the structure (or Reeb) vector field on $M$. If $X$ is a vector field on $M$ we write $JX=\phi X+\eta(X)N$, where $\phi X$ denotes the tangent component of $JX$. Then $\eta(X)=g(X,\xi)$, $\phi$ is called the structure tensor on $M$ and $(\phi , \xi , \eta , g)$ is an almost contact metric structure on $M$ induced by the Kaehlerian structure of $\C P^m$. The classification of homogeneous real hypersurfaces in $\C P^m$ was obtained by Takagi, see \cite{K1}, \cite{T1}, \cite{T2}, \cite{T3}. His classification contains 6 types of real hypersurfaces. Among them we find type $(A_1)$ real hypersurfaces that are geodesic hyperspheres of radius $r$, $0 < r < \frac {\pi}{2}$, and type $(A_2)$ real hypersurfaces that are tubes of radius $r$, $0 < r <\frac {\pi}{2}$, over totally geodesic complex projective spaces $\C P^n$, $0 < n < m-1$. We will call both types of real hypersurfaces type $(A)$ real hypersurfaces. They are Hopf, that is, the structure vector field is principal, and are the unique real hypersurfaces in $\C P^m$ such that $A\phi =\phi A$, see \cite{O}. 

Ruled real hypersurfaces in $\C P^m$ satisfy that the maximal holomorphic distribution on $M$, $\D$, given at any point by the vectors orthogonal to $\xi$, is integrable and its integral manifolds are totally geodesic $\C P^{m-1}$. Equivalently, $g(A\D ,\D)=0$. For examples of ruled real hypersurfaces see \cite{K2} or \cite{LR}.

The Tanaka-Webster connection, \cite{TA}, \cite{W}, is the canonical affine connection defined on a non-degenerate, pseudo-Hermitian CR-manifold. As a generalization of this connection, Tanno, \cite{TAN}, defined the generalized Tanaka-Webster connection for contact metric manifolds by

\begin{multi} \label{1.1}
\hat{\nabla}_XY=\nabla_XY+(\nabla_X \eta)(Y)\xi -\eta(Y)\nabla_X \xi -\eta(X)\phi Y
\end{multi}

\noindent for any vector fields $X,Y$ on the manifold.

Using the almost contact metric structure on $M$ and the naturally extended affine connection of Tanno's generalized Tanaka-Webster connection, Cho defined the $k$-th generalized Tanaka-Webster connection $\hat{\nabla}^{(k)}$ for a real hypersurface $M$ in $\C P^m$, see \cite{CH1}, \cite{CH2}, by

\begin{multi} \label{1.2}
\hat{\nabla}^{(k)}_XY=\nabla_XY+g(\phi AX,Y)\xi -\eta (Y)\phi AX-k\eta (X)\phi Y
\end{multi}

\noindent for any X,Y tangent to $M$ where $k$ is a non-zero real number. Then $\hat{\nabla}^{(k)}\eta =0$, $\hat{\nabla}^{(k)}\xi =0$, $\hat{\nabla}^{(k)}g=0$, $\hat{\nabla}^{(k)}\phi =0$. In particular, if the shape operator of a real hypersurface satisfies $\phi A+A \phi =2k\phi$, the $k$-th generalized Tanaka-Webster connection coincides with the Tanaka-Webster connection.

Here we can consider the tensor field of type (1,2) given by the difference of both connections $F^{(k)}(X,Y)=g(\phi AX,Y)\xi -\eta(Y)\phi AX-k\eta(X)\phi Y$, for any $X,Y$ tangent to $M$, see \cite{KN} Proposition 7.10, pages 234-235. We will call this tensor the $k$-th Cho tensor on $M$. Associated to it, for any $X$ tangent to $M$ and any nonnull real number $k$ we can consider the tensor field of type (1,1) $F_X^{(k)}$, given by $F_X^{(k)}Y=F^{(k)}(X,Y)$ for any $Y \in TM$. This operator will be called the $k$-th Cho operator corresponding to $X$. Notice that if $X \in \D$, the corresponding Cho operator does not depend on $k$ and we simply write it $F_X$. The torsion of the connection $\hat{\nabla}^{(k)}$ is given by $T^{(k)}(X,Y)=F_X^{(k)}Y-F_Y^{(k)}X$ for any $X,Y$ tangent to $M$. We define the $k$-th torsion operator associated to $X$ to the operator given by $T_X^{(k)}Y=T^{(k)}(X,Y)$, for any $X,Y$ tangent to $M$.

Let $\cal{L}$ denote the Lie derivative on $M$. Therefore, ${\cal L}_XY=\nabla_XY-\nabla_YX$ for any $X,Y$ tangent to $M$. Now we can define on $M$ a differential operator of first order, associated to the $k$-th generalized Tanaka-Webster connection, given by

\[  {\cal L}_X^{(k)}Y=\hat{\nabla}_X^{(k)}Y-\hat{\nabla}_Y^{(k)}X={\cal L}_XY+T_X^{(k)}Y       \]

\noindent for any $X,Y$ tangent to $M$. We will call it the derivative of Lie type associated to the $k$-th generalized Tanaka-Webster connection.

Let now $L$ be a symmetric tensor of type (1,1) defined on $M$. We can consider then the type (1,2) tensor $L_F^{(k)}$ associated to $L$ in the following way

\[         L_F^{(k)}(X,Y)=\lbrack F_X^{(k)},L \rbrack Y=F^{(k)}_XLY-LF^{(k)}_XY                  \]

\noindent for any $X,Y$ tangent to $M$. We also can consider another tensor of type (1,2), $L_T^{(k)}$, associated to $L$, by

\[       L_T^{(k)}(X,Y)=\lbrack T_X^{(k)},L \rbrack Y=T_X^{(k)}LY-LT_X^{(k)}Y              \]

\noindent for any $X,Y$ tangent to $M$. Notice that if $X \in  \D$, $L_F^{(k)}$ does not depend on $k$. We will write it simply $L_F$.
 
In \cite{PSU}, respectively, \cite{P} we proved non-existence of real hypersurfaces in $\C P^m$, $m \geq 3$, such that the tensors of type (1,2) associated to the shape operator, $A_F^{(k)}=0$, respectively, $A_T^{(k)}=0$, for any nonnull real number $k$. Further results on such tensors were obtained in \cite{PPL} and \cite{PPL2}.

The purpose of the present paper is to study the behaviour of both tensors with respect to the structure operator $\phi$. We will say that $A_F^{(k)}$ is pure with respect to $\phi$ if $A_F^{(k)}(\phi X,Y)=A_F^{(k)}(X,\phi Y)$, for any $X,Y$ tangent to $M$, see \cite{SA}, \cite{SAA}. We will say that $A_F$ is $\eta$-pure with respect to $\phi$ if $A_F(\phi X,Y)=A_F(X,\phi Y)$, for any $X,Y \in \D$. Analogously we will say that $A_F$ is $\eta$-skewpure with respect to $\phi$ if $A_F(\phi X,Y)+A_F(X,\phi Y)=0$. for any $X,Y \in \D$. We will prove the following

\noindent {\bf Theorem 1} {\it Let $M$ be a real hypersurface in $\C P^m$, $m\geq 3$. Then $A_F$ is $\eta$-pure with respect to $\phi$ if and only if $M$ is locally congruent to a ruled real hypersurface.}

Also we will prove

\noindent {\bf Theorem 2} {\it Let $M$ be a real hypersurface in $\C P^m$, $m \geq 3$. Thn $A_F$ is $\eta$-skewpure with respect to $\phi$ if and only if $M$ is locally congruent to one of the following real hypersurfaces:

1) A tube of radius $\frac{\pi}{4}$ around a complex submanifold of $\C P^m$.

2) A real hypersurface of type $(A)$.

3) A ruled real hypersurface.}

On the other hand we also have

\noindent {\bf Theorem 3} {\it Let $M$ be a real hypersurface in $\C P^m$, $m \geq 3$. Then $A_F(\phi X,Y)=\phi A_F(X,Y)$ for any $X,Y \in \D$ if and only if $M$ is locally congruent to a ruled real hypersurface.}

Concerning the tensor $A_T^{(k)}$ we will prove

\noindent {\bf Theorem 4} {\it There does not exist any real hypersurface in $\C P^m$, $m \geq 3$, such that $A_T^{(k)}$ is $\eta$-pure with respect to $\phi$, for any nonnull real number $k$.}

Also we will obtain

\noindent {\bf Theorem 5} {\it Let $M$ be a real hypersurface in $\C P^m$, $m \geq 3$, and $k$ a nonnull real number. Then $A_T^{(k)}$ is $\eta$-skewpure with respect to $\phi$ if and only if $M$ is locally congruent to a real hypersurface of type $(A)$.}

As in the case of $A_F$ we can prove

\noindent {\bf Theorem 6} {\it Let $M$ be a real hypersurface in $\C P^m$, $m \geq 3$, and $k$ a nonnull real number. Then $A_T^{(k)}(\phi X,Y)=\phi A_T^{(k)}(X,Y)$, for any $X,Y \in \D$, if and only if $M$ is locally congruent to a ruled real hypersurface.}

And

\noindent {\bf Theorem 7} {\it There does not exist any real hypersurface in $\C P^m$, $m \geq 3$, such that $A_T^{(k)}(X,\phi Y)=\phi A_T^{(k)}(X,Y)$ for any $X,Y \in \D$ and any nonnull real number $k$.}

{\bf Aknowledgements} This work was supported by MICINN Project PID 2020-116126GB-100 and Project PY20-01391 from Junta de Andaluc\'{i}a.

\section{Preliminaries.}

Throughout this paper, all manifolds, vector fields, etc., will be considered of class $C^{\infty}$ unless otherwise stated. Let $M$ be a connected real hypersurface in $\C P^m$, $m\geq 2$, without boundary. Let $N$ be a locally defined unit normal vector field on $M$. Let $\nabla$ be the Levi-Civita connection on $M$ and $(J,g)$ the Kaehlerian structure of $\C P^m$.

For any vector field $X$ tangent to $M$ we write $JX=\phi X+\eta(X)N$, and $-JN=\xi$. Then $(\phi,\xi,\eta,g)$ is an almost contact metric structure on $M$, see \cite{B}. That is, we have

\begin{equation} \label{2.1}
\phi^2X=-X+\eta(X)\xi,  \quad  \eta(\xi)=1,  \quad   g(\phi X,\phi Y)=g(X,Y)-\eta(X)\eta(Y)
\end{equation}

\noindent for any vectors $X,Y$ tangent to $M$. From (\ref{2.1}) we obtain

\begin{multi} \label{elem}
\phi\xi =0,   \quad   \eta(X)=g(X,\xi).
\end{multi}

\noindent From the parallelism of $J$ we get

\begin{equation} \label{lasphi}
(\nabla_X\phi)Y=\eta(Y)AX-g(AX,Y)\xi
\end{equation}

\noindent and

\begin{equation} \label{elem2}
\nabla_X\xi =\phi AX
\end{equation}

\noindent for any $X,Y$ tangent to $M$, where $A$ denotes the shape operator of the immersion. As the ambient space has holomorphic sectional curvature $4$, the equation Codazzi is given by

\begin{equation} \label{Codazzi}
(\nabla_XA)Y-(\nabla_YA)X=\eta(X)\phi Y-\eta(Y)\phi X-2g(\phi X,Y)\xi
\end{equation}

\noindent for any vectors $X,Y$ tangent to $M$. We will call the maximal holomorphic distribution $\D$ on $M$ to the following one: at any $p \in M$, $\D (p)=\{ X\in T_pM \vert  g(X,\xi)=0\}$. We will say that $M$ is Hopf if $\xi$ is principal, that is, $A\xi =\alpha\xi$ for a certain function $\alpha$ on $M$.

In the sequel we need the following result:

\noindent {\bf Theorem 2.1, \cite{M}.} {\it
If $\xi$ is a principal curvature vector with corresponding principal curvature $\alpha$ and $X \in \D$ is principal with principal curvature $\lambda$, then $2\lambda -\alpha \neq 0$ and $\phi X$ is principal with principal curvature $\frac {\alpha\lambda +2}{2\lambda -\alpha}$.}

\section{Proofs of results concerning $A_F$.}

In order to prove Theorem 1, we should have $F_{\phi X}AY-AF_{\phi X}Y=F_XA\phi Y-AF_X\phi Y$, for any $X,Y \in \D$. This yields

\begin{multi} \label{3.1}
g(\phi A\phi X,AY)\xi -\eta(AY)\phi A\phi X-g(\phi A\phi X,Y)A\xi     \\                                                                                                     =g(\phi AX,A\phi Y)\xi -\eta(A\phi Y)\phi AX-g(AX,Y)A\xi        \\
\end{multi}

\noindent for any $X,Y \in \D$. If $M$ is Hopf with $A\xi =\alpha\xi$, the scalar product of (\ref{3.1}) and $\xi$ gives $g(\phi A\phi X,AY)-\alpha g(\phi A\phi X,Y)=g(\phi AX,A\phi Y)-\alpha g(\phi AX,Y)$ for any $X,Y \in \D$. Let us suppose that $X \in \D$ satisfies $AX=\lambda X$. Then $A\phi X=\mu\phi X$, $\mu =\frac {\alpha\lambda +2}{2\lambda -\alpha}$, and we obtain $-\lambda\mu +\alpha\mu =\lambda\mu -\alpha\lambda$. That is, $2\lambda\mu =\alpha(\mu +\lambda)$. This implies $\frac{2\alpha\lambda^2+4\lambda}{2\lambda -\alpha}=\alpha\left(\frac{\alpha\lambda +2}{2\lambda -\alpha} +\lambda\right)=\alpha\left(\frac{2(1+\lambda)^2}{2\lambda -\alpha}\right)$. Thus $\alpha\lambda^2+2\lambda =\alpha\lambda^2+\alpha$, and so, $\lambda =\frac{\alpha}{2}$. As $2\lambda\mu =\alpha(\mu +\lambda)$ we get $\alpha\mu =\alpha(\mu +\lambda)$. Then $\alpha\lambda =\frac{\alpha^2}{2}=0$, that is, $\alpha =0$ and also $\lambda =0$, a contradiction with the fact of $2\lambda -\alpha \neq 0$.

This means that $M$ must be non Hopf. Therefore, locally we can write $A\xi =\alpha\xi +\beta U$, $U$ being a unit vector field in $\D$, $\alpha$ and $\beta$ functions on $M$ and $\beta \neq 0$. We also define ${\D}_U$ as the orthogonal complementary distribution in $\D$ to the one spanned by $U$ and $\phi U$. With this in mind (\ref{3.1}) becomes

\begin{multi} \label{3.2}
g(\phi A\phi X,AY)\xi -\beta g(Y,U)\phi A\phi X-g(\phi A\phi X,Y)A\xi       \\                                                                                              =g(\phi AX,A\phi Y)\xi -\beta g(\phi Y,U)\phi AX-g(AX,Y)A\xi    \\
\end{multi}

\noindent for any $X,Y \in \D$. The scalar product of (\ref{3.2}) and $\phi U$ gives $-\beta g(Y,U)g(A\phi X,U)=-\beta g(\phi Y,U)g(AX,U)$, for any $X,Y \in \D$. Taking $Y=U$ we obtain $-\beta g(AU,\phi X)=0$ for any $X \in \D$. As we suppose $\beta \neq 0$ and changing $X$ by $\phi X$ we have $g(AU,X)=0$ for any $X \in \D$. This means that
 
\begin{multi} \label{3.3}
AU=\beta\xi.
\end{multi}

The scalar product of (\ref{3.2}) and $U$ yields $-\beta g(Y,U)g(\phi A\phi X,U)-\beta g(\phi A\phi X,Y)=-\beta g(\phi Y,U)g(\phi AX,U)-\beta g(AX,Y)$, for any $X,Y \in \D$. As $\beta \neq 0$ we have

\begin{multi} \label{3.4}
g(Y,U)g(A\phi U,\phi X)+g(A\phi Y,\phi X)=g(\phi Y,U)g(A\phi U,X)-g(AX,Y)
\end{multi}

\noindent for any $X,Y \in \D$. If we take $X=U$ in (\ref{3.4}) it follows $2g(A\phi U,\phi X)=-g(AU,X)$ for any $X \in \D$. From (\ref{3.3}), changing $X$ by $\phi X$, we obtain $g(A\phi U,X)=0$ for any $X \in \D$. Therefore

\begin{multi} \label{3.5}
A\phi U=0.
\end{multi}

Now the scalar product of (\ref{3.2}) and $Z \in {\D}_U$ implies $-\beta g(Y,U)g(\phi A\phi X,Z)=-\beta g(\phi Y,U)g(\phi AX,Z)$, for any $X,Y \in \D$, $Z \in {\D}_U$. If $Y=\phi U$ we obtain $\beta g(\phi AX,Z)$ $=0$ for any $X \in \D$, $Z \in {\D}_U$. If we change $Z$ by $\phi Z$ and bear in mind that $\beta \neq 0$, it follows $g(AZ,X)=0$ for any $Z \in {\D}_U$, $X \in \D$. Therefore,

\begin{multi} \label{3.6}
AZ=0
\end{multi}

\noindent for any $Z \in {\D}_U$. From (\ref{3.3}), (\ref{3.5}) and (\ref{3.6}), $M$ is locally congruent to a ruled real hypersurface. The converse is trivial and we have finished the proof of Theorem 1.

Now if $A_F$ is $\eta$-skewpure we have

\begin{multi} \label{3.7}
g(\phi A\phi X,AY)\xi -\eta(AY)\phi A\phi X-g(\phi A\phi X)A\xi +g(\phi AX,A\phi Y)\xi         \\                                                                             -\eta(A\phi Y)\phi AX-g(AX,Y)A\xi =0      \\
\end{multi}

\noindent for any $X,Y \in \D$. Let us suppose that $M$ is Hopf and write $A\xi =\alpha\xi$. If we take the scalar product of (\ref{3.7}) and $\xi$ it follows $g(\phi A\phi X,AY)-\alpha g(\phi A\phi X,Y)+g(\phi AX,A\phi Y)-\alpha g(AX,Y)=0$, for any $X,Y \in \D$. This means that $A\phi A\phi X-\alpha\phi A\phi X-\phi A\phi AX-\alpha AX=0$ for any $X \in \D$. If we take $X \in \D$ such that $AX=\lambda X$, as $A\phi X=\mu\phi X$, we get $-\lambda\mu +\alpha\mu +\lambda\mu -\alpha\lambda =0$. That is, $\alpha(\mu -\lambda)=0$. Thus either $\alpha =0$, and by \cite{CR}   we have 1) in Theorem 2, or $\mu =\lambda$. This means that $A\phi =\phi A$ and in this case we have 2) in Theorem 2.

If $M$ is non Hopf following the same steps as in Theorem 1 we obtain 3) in Theorem 2, finishing its proof.

If we suppose that $M$ satisfies the condition in Theorem 3 we must have $F_{\phi X}AY-AF_{\phi X}Y=\phi F_XAY-\phi AF_XY$, for any $X,Y \in \D$. This yields

\begin{multi} \label{3.8}
g(\phi A\phi X,AY)\xi -\eta(AY)\phi A\phi X-g(\phi A\phi X,Y)A\xi =-\eta(AY)\phi^2AX          \\                                                                             -g(\phi AX,Y)\phi A\xi           \\
\end{multi}

\noindent for any $X,Y \in \D$. If we suppose that $M$ is Hopf, the scalar product of (\ref{3.8}) and $\xi$ gives $g(\phi A\phi X,AY)-\alpha g(\phi A\phi X,Y)=0$. Therefore, $A\phi A\phi X-\alpha\phi A\phi X=0$, for any $X \in \D$. If we suppose that $X \in \D$ satisfies $AX=\lambda X$ we obtain $\mu (\alpha -\lambda)=0$. Therefore, either $\mu =0$ and then $\alpha \neq 0$ and $\lambda =-\frac{2}{\alpha}$, or if $\mu \neq 0$, $\alpha =\lambda$ and then $\mu =\frac{\alpha^2+2}{\alpha}$. Moreover all principal curvatures are constant and by \cite{K1}, $M$ must be locally congruent to a real hypersurface appearing among the six types in Takagi's list. Looking at such types, none has our principal curvatures, \cite{T2}, proving that $M$ must be non Hopf.

We write as above $A\xi =\alpha\xi +\beta U$, with the same conditions. Then (\ref{3.8}) becomes

\begin{multi} \label{3.9}
g(A\phi A\phi X,Y)\xi -\beta g(Y,U)\phi A\phi X-g(\phi A\phi X,Y)A\xi                    \\                                                                                        =-\beta g(Y,U)\phi^2AX-\beta g(\phi AX,Y)\phi U         \\
\end{multi}

\noindent for any $X,Y \in \D$. The scalar product of (\ref{3.9}) and $\phi U$ gives, bearing in mind that $\beta \neq 0$,

\begin{multi} \label{3.10}
g(Y,U)g(AU,\phi X)=g(Y,U)g(\phi AX,U)-g(\phi AX,Y),                                                                                                                                       
\end{multi}

\noindent for any $X,Y \in \D$. If $X=Y=U$, we get $g(AU,\phi U)=-2g(AU,\phi U)$. Thus

\begin{multi} \label{3.11}
g(AU,\phi U)=0.
\end{multi}

If we take $Y=U$, $X \in \D$ and orthogonal to $U$ in (\ref{3.10}), we have $g(\phi AU,X)=2g(A\phi U,X)$, for such an $X$. From (\ref{3.11}) the same is true for $X=U$. Therefore, $2A\phi U-\phi AU$ has no component in $\D$. As its scalar product with $\xi$ also vanishes we get

\begin{multi} \label{3.12}
2A\phi U=\phi AU.
\end{multi}

If we take $Y=\phi U$, $X \in \D$ in (\ref{3.10}) it follows $g(AX,U)=0$ for any $X \in \D$. Thus

\begin{multi} \label{3.13}
AU =\beta\xi
\end{multi}

\noindent and from (\ref{3.11}),

\begin{multi} \label{3.14}
A\phi U=0.
\end{multi}

The scalar product of (\ref{3.9}) and $U$, bearing in mind (\ref{3.13}) and (\ref{3.14}), gives $g(\phi A\phi X,$ $Y)=g(Y,U)g(\phi^2AX,U)=-g(Y,U)g(AX,U)=0$, for any $X \in \D$. Taking $\phi X \in {\D}_U$ instead of $X$ we obtain $\phi AX=0$. Applying $\phi$ we get

\begin{multi} \label{3.15}
AX=0
\end{multi}

\noindent for any $X \in {\D}_U$. From (\ref{3.13}), (\ref{3.14}) and (\ref{3.15}), $M$ must be locally congruent to a ruled real hypersurface and we have finished the proof of Theorem 3.

{\bf Remark 1.-} With similar proofs to the proof of Theorem 3, we can obtain other characterizations of ruled real hypersurfaces in $\C P^m$, $m \geq 3$, if we consider any of the following conditions:

1) $A_F(\phi X,Y)+\phi A_F(X,Y)=0$, for any $X,Y \in \D$.

2) $A_F(X,\phi Y)=\phi A_F(X,Y)$, for any $X, Y \in \D$.

3) $A_F(X,\phi Y)+\phi A_F(X,Y)=0$, for any $X,Y \in \D$.

\section{Results concerning $A_T^{(k)}$.}

If we suppose that $A_T^{(k)}$ is $\eta$-pure with respect to $\phi$ we will have $F_{\phi X}AY-F_{AY}^{(k)}\phi X-AF_{\phi X}Y+AF_Y\phi X=F_XA\phi Y-F_{A\phi Y}^{(k)}X-AF_X\phi Y+AF_{\phi Y}X$, for any $X,Y \in \D$. This yields

\begin{multi} \label{4.1}
g(\phi A\phi X ,AY)\xi -\eta(AY)\phi A\phi X -g(A^2Y,X)\xi -k\eta(AY)X-g(\phi A\phi X ,Y)A\xi             \\                                                    +g(AY,X)A\xi =g(\phi AX,A\phi Y)\xi -\eta(A\phi Y)\phi AX-g(\phi A^2\phi Y,X)\xi       \\                                                                             +k\eta(A\phi Y)\phi X -g(AX,Y)A\xi +g(\phi A\phi X)A\xi  \\
\end{multi}

\noindent for any $X,Y \in \D$. Let us suppose that $M$ is Hopf with $A\xi =\alpha\xi$. Then (\ref{4.1}) becomes $g(\phi A\phi X,Y)\xi -g(A^2Y,X)\xi -\alpha g(\phi A\phi X,Y)\xi +\alpha g(AY,X)\xi =g(\phi AX,A\phi Y)\xi -g(\phi A^2\phi Y,X)\xi -\alpha g(AX,Y)\xi +\alpha g(\phi A\phi Y,X)\xi$, for any $X,Y \in \D$. Let us suppose that $X \in \D$ satisfies $AX=\lambda X$. Then $A\phi X=\mu\phi X$ and from last equation we obtain $-\lambda\mu-\lambda^2+2\alpha\mu +2\alpha\lambda =\lambda\mu +\mu^2$. That is $(\mu +\lambda)^2-2\alpha(\mu +\lambda)=0$. Thus $(\mu +\lambda)(\mu +\lambda -2\alpha)=0$. If $\mu +\lambda =0$, as $\mu =\frac{\alpha\lambda +2}{2\lambda -\alpha}$, we get $2\lambda^2+2=0$, which is impossible. Therefore $\mu +\lambda =2\alpha$ and the value of $\mu$ yields $\lambda^2-2\alpha\lambda +1+\alpha^2=0$. This equation has not real solutions and this implies that our real hypersurface must be non Hopf.

As in the previous section we write locally $A\xi =\alpha\xi +\beta U$, with the same conditions, and also make the following computations locally. The scalar product of (\ref{4.1}) and $\phi U$ gives $-\eta(AY)g(A\phi X,U)-k\eta(AY)g(X,\phi U)=-\eta(A\phi Y)g(AX,U)+k\eta(A\phi Y)g(X,U)$, for any $X,Y \in \D$. That is, bearing in mind that $\beta \neq 0$,

\begin{multi} \label{4.2}
g(Y,U)g(A\phi X,U)+kg(Y,U)g(X,\phi U)=g(\phi Y,U)g(AX,U)-kg(\phi Y,U)g(X,U)
\end{multi}

\noindent for any $X,Y \in \D$. Take $Y=\phi U$ in (\ref{4.2}) to obtain $g(AX,U)-kg(X,U)=0$, for any $X \in \D$. Therefore,

\begin{multi} \label{4.3}
AU=\beta\xi +kU.
\end{multi}

Now the scalar product of (\ref{4.1}) and $U$ yields

\begin{multi} \label{4.4}
-g(Y,U)g(\phi A\phi X,U)-kg(y,U)g(X,U)-g(\phi A\phi X,Y)+g(AY,X)  \\
=-g(\phi Y,U)g(\phi AX,U)+kg(\phi Y,U)g(\phi X,U)-g(AX,Y)+g(\phi A\phi Y,X)  \\
\end{multi}

\noindent for any $X,Y \in \D$. Taking $Y=U$ in (\ref{4.4}) we obtain

\begin{multi} \label{4.5}
-kg(X,U)-3g(\phi A\phi X,U)+2g(AU,X)=0
\end{multi}

\noindent for any $X \in \D$. Taking $X \in {\D}_U$ and changing $X$ by $\phi X$ in (\ref{4.5}) we get $g(A\phi U,X)=0$, for any $X \in {\D}_U$. If $X=U$ in (\ref{4.5}) we have $-k+3g(A\phi U,\phi U)+2k =0$. Bearing in mind (\ref{4.3}) we have obtained

\begin{multi} \label{4.6}
A\phi U=-\frac{k}{3} \phi U.
\end{multi}

Moreover, the scalar product of (\ref{4.1}) and $\phi Z \in {\D}_U$ implies

\begin{multi} \label{4.7}
-g(Y,U)g(A\phi X,Z)-kg(Y,U)g(X,\phi Z)=-g(\phi Y,U)g(AX,Z)+kg(\phi Y,U)g(X,Z)
\end{multi}

\noindent for any $X,Y \in \D$, $Z \in {\D}_U$. Taking $Y=\phi U$ in (\ref{4.7}) we obtain $g(AZ,X)-kg(Z,X)=0$, for any $Z \in {\D}_U$, $X \in \D$, and this yields

\begin{multi} \label{4.8}
AZ=kZ
\end{multi}

\noindent for any $Z \in {\D}_U$. Take $Z \in {\D}_U$. Then $AZ=kZ$ and $A\phi Z=k\phi Z$. From Codazzi equation, $\nabla_{\phi Z}(kZ)-A\nabla_{\phi Z}Z-\nabla_Z(k\phi Z)+A\nabla_Z\phi Z=2\xi$. Its scalar product with $\xi$  gives $-kg(Z,\phi A\phi Z)-g(\nabla_{\phi Z}Z,\alpha\xi +\beta U)+kg(\phi Z,\phi AZ)+g(\nabla_Z\phi Z,\alpha\xi +\beta U)=2$. Then, $\beta g(\lbrack Z,\phi X \rbrack ,U)+2k^2+\alpha g(Z,\phi A\phi Z)-\alpha g(\phi Z,\phi AZ)=2$. Therefore

\begin{multi} \label{4.9}
g(\lbrack Z, \phi Z \rbrack ,U)=\frac{2-2k^2+2\alpha k}{\beta}.
\end{multi}

And its scalar product with $U$ implies $-kg(\lbrack Z,\phi Z \rbrack ,U)-g(\nabla_{\phi Z}Z,\beta\xi +kU)+g(\nabla_Z\phi Z,\beta\xi +kU)=0$. This gives $\beta g(Z,\phi A\phi Z)-\beta g(\phi Z,\phi AZ)=0$ or $2\beta k=0$, which is impossible and proves Theorem 4.

Suppose now that $A_T^{(k)}$ is $\eta$-skewpure with respect to $\phi$. We should have

\begin{multi} \label{4.10}
g(\phi A\phi X,AY)\xi -\eta(AY)\phi A\phi X-g(A^2Y,X)\xi -k\eta(AY)X+g(\phi AX,A\phi Y)\xi         \\
-\eta(A\phi Y)\phi AX-g(\phi A^2\phi Y,X)\xi +k\eta(A\phi Y)\phi X=0
\end{multi}

\noindent for any $X,Y \in \D$. Let us suppose that $M$ is Hopf. Then (\ref{4.10}) gives $g(A\phi A\phi X,Y)-g(A^2X,Y)-g(\phi A\phi AX,Y)-g(\phi A^2\phi X,Y)=0$, for any $X,Y \in \D$. Then $A\phi A\phi X-A^2X-\phi A\phi AX-\phi A^2\phi X =0$, for any $X \in \D$. If $X \in \D$ satisfies $AX=\lambda X$ we obtain $-\lambda\mu -\lambda^2+\lambda\mu +\mu^2=0$. Therefore, $\lambda^2=\mu^2$. As in the previous Theorem $\lambda +\mu =0$ gives a contradiction. This means that $\lambda =\mu$ and $\phi A=A\phi$. This yields that $M$ must be locally congruent to a real hypersurface of type $(A)$. The converse is immediate.

Suppose then that $M$ is non Hopf and $A\xi =\alpha\xi +\beta U$. Taking the scalar product of (\ref{4.10}) and $\phi U$ we have

\begin{multi} \label{4.11}
-g(Y,U)g(A\phi X,U)-kg(Y,U)g(X,\phi U)-g(\phi Y,U)g(AX,U)   \\                                                                                                               +kg(\phi Y,U)g(X,U)=0                  \\
\end{multi}

\noindent for any $X,Y \in \D$. If $Y=\phi U$ in (\ref{4.11}) we obtain $g(AU,X)-kg(U,X)=0$ for any $X \in \D$ and this yields

\begin{multi} \label{4.12}
AU=\beta\xi +kU.
\end{multi}

Following the above proof step by step we can also see that $A\phi U=k\phi U$ and $AZ=kZ$, for any $Z \in {\D}_U$. If we apply again Codazzi equation to $Z$ and $\phi Z$, $Z \in {\D}_U$ we obtain $k\beta =0$, which is impossible and finishes the proof of Theorem 5.

Condition in Theorem 6 implies

\begin{multi} \label{4.13}
g(\phi A\phi X,AY)\xi -\eta(AY)\phi A\phi X-g(A^2Y,X)\xi -g(A^2Y,X)\xi   \\                                                                                                        -g(\phi A\phi X,Y)A\xi +g(AX,Y)A\xi =-\eta(AY)\phi^2AX-g(\phi AX,Y)\phi A\xi      \\                                                                                            +g(\phi AY,X)\phi A\xi         \\
\end{multi}

\noindent for any $X,Y \in \D$. If $M$ is Hopf, (\ref{4.13}) yields $g(A\phi A\phi X,Y)\xi -g(A^2X,Y)\xi -\alpha g(\phi A\phi X,Y)\xi -\alpha g(AX,Y)\xi =0$, for any $X,Y \in \D$. Its scalar product with $\xi$ shows that $A\phi A\phi X-A^2X-\alpha\phi A\phi X+\alpha AX=0$, for any $X \in \D$: If $X \in \D$ satisfies $AX=\lambda X$, $A\phi X=\mu\phi X$ and we obtain $(\alpha -\lambda)(\mu +\lambda)=0$. As we saw before, $\lambda +\mu \neq 0$. Therefore, $\lambda =\alpha$ and as $\mu =\frac{\alpha\lambda +2}{2\lambda -\alpha}$, we get $\mu =\frac{\alpha^2+2}{\alpha}$. Thus $M$ has two distinct constant principal curvatures. From \cite{CR} $M$ must be locally congruent to a geodesic hypersphere. In this case $M$ has only a principal curvature on $\D$. That means that $\alpha =\frac{\alpha^2+2}{\alpha}$, which is impossible. Therefore $M$ must be non Hopf and as above, we write $A\xi =\alpha\xi +\beta U$. In this case (\ref{4.13}) looks like 

\begin{multi} \label{4.14}
g(\phi A\phi X,AY)\xi -\beta g(Y,U)\phi A\phi X)-g(A^2Y,X)\xi-g(\phi A\phi X)A\xi        \\                                                                                  +g(AX,Y)A\xi =-\beta g(Y,U)\phi^2AX-\beta g(\phi AX,Y)\phi U+\beta g(\phi AY,X)\phi U   \\
\end{multi}

\noindent for any $X,Y \in \D$. Bearing in mind that $\beta \neq 0$, the scalar product of (\ref{4.14}) and $\phi U$ gives $-g(Y,U)g(A\phi X,U)=-g(Y,U)g(\phi AX,U)-g(\phi AX,Y)+g(\phi AY,X)$, for any $X,Y \in \D$. Taking $Y=U$ we get $g(\phi AU,X)=-g(\phi AX,U)-g(\phi AX,U)+g(\phi AU,X)$. Therefore, $g(A\phi U,X)=0$ for any $X \in \D$, and

\begin{multi} \label{4.15}
A\phi U=0.
\end{multi}

Taking $Y=\phi U$ in the last equation we get $0=-g(AX,U)+g(\phi A\phi U,X)$, for any $X \in \D$. Bearing in mind (\ref{4.15}) this implies $g(AU,X)=0$ for any $X \in \D$ and so

\begin{multi} \label{4.16}
AU=\beta\xi.
\end{multi}

From (\ref{4.15}) and (\ref{4.16}), ${\D}_U$ is $A$-invariant. If in the equality used to find (\ref{4.15}) and (\ref{4.16}) we take $Y \in {\D}_U$ we obtain $0=g(\phi AY,X)+g(A\phi Y,X)$, for any $Y \in {\D}_U$, $X \in \D$. Then, $\phi AY+A\phi Y=0$ for any $Y \in {\D}_U$. If we suppose that $AY=\lambda Y$, we get $A\phi Y=-\lambda\phi Y$.

The scalar product of (\ref{4.14}) and $Z \in {\D}_U$ gives

\begin{multi} \label{4.17}
-g(Y,U)g(\phi A\phi X,Z)=g(Y,U)g(AX,Z)
\end{multi}

\noindent for any $X,Y \in \D$, $Z \in {\D}_U$. Taking $Y=U$, $X=Z$ we obtain $g(A\phi Z,\phi Z)=g(AZ,Z)$. This yields $\lambda=-\lambda$. Therefore $\lambda =0$ and $M$ is locally congruent to a ruled real hypersurface. This finishes the proof of Theorem 6.

The condition in Theorem 7 yields

\begin{multi} \label{4.18}
g(\phi AX,A\phi Y)\xi -\eta(A\phi Y)\phi AX-g(\phi A^2\phi Y,X)\xi +k\eta(A\phi Y)\phi X   \\
-g(AX,Y)A\xi +g(\phi A\phi Y,X)A\xi =-\eta(AY)\phi^2AX-k\eta(AY)X    \\
-g(\phi AX,Y)\phi A\xi +g(\phi AY,X)\phi A\xi  \\
\end{multi}

\noindent for any $X,Y \in \D$. If we suppose that $M$ is Hopf with $A\xi =\alpha\xi$, (\ref{4.18}) becomes

\begin{multi} \label{4.19}
g(\phi AX,A\phi Y)\xi -g(\phi A^2\phi Y,X)\xi -\alpha g(AX,Y)\xi +\alpha g(\phi A\phi Y,X)\xi =0
\end{multi}

\noindent for any $X,Y \in \D$. This gives $-\phi A\phi AX-\phi A^2\phi X-\alpha AX+\alpha\phi A\phi X=0$ for any $X \in \D$. If $X \in \D$ satisfies $AX=\lambda X$ we obtain $(\mu -\alpha)(\lambda + \mu)=0$. As in previous Theorems this case leads to a contradiction.

Thus $M$ must be non Hopf and, as usual, we write $A\xi =\alpha\xi +\beta U$. In this case (\ref{4.18}) implies

\begin{multi} \label{4.20}
g(\phi AX,A\phi Y)\xi -\beta g(\phi Y,U)\phi AX-g(\phi A^2\phi Y,X)\xi +k\beta g(\phi Y,U)\phi X    \\
-g(AX,Y)A\xi +g(\phi A\phi Y,X)A\xi =\beta g(Y,U)AX- \beta^2g(Y,U)g(X,U)\xi      \\                                                                                                   -k\beta g(Y,U)X-\beta g(\phi AX,Y)\phi U+\beta g(\phi AY,X)\phi U    \\
\end{multi}

\noindent for any $X,Y \in \D$. The scalar product of (\ref{4.20}) and $\phi U$, bearing in mind that $\beta \neq 0$ gives

\begin{multi} \label{4.21}
-g(\phi Y,U)g(AX,U)+kg(\phi Y,U)g(X,U)=g(Y,U)g(A\phi U,X)     \\
-kg(Y,U)g(\phi U,X)-g(\phi AX,Y)+g(\phi AY,X)    \\
\end{multi}

\noindent for any $X,Y \in \D$. If we take $Y \in {\D}_U$ in (\ref{4.21}) we obtain

\begin{multi} \label{4.22}
g(AX,\phi Y)-g(A\phi X,Y)=0
\end{multi}

\noindent for any $X \in \D$, $Y \in {\D}_U$.

Taking $X \in {\D}_U$ in (\ref{4.21}) it follows

\begin{multi} \label{4.23}
-g(\phi Y,U)g(AU,X)=g(Y,U)g(A\phi U,X)+g(A\phi Y,X)+g(\phi AY,X)
\end{multi}

\noindent for any $X \in {\D}_U$, $Y \in \D$.

Take $X=U$ and $\phi Y$ instead of $Y$ in (\ref{4.22}) to have

\begin{multi} \label{4.24}
g(AU,Y)+g(A\phi U,\phi Y)=0
\end{multi}

\noindent for any $Y \in {\D}_U$. Take $Y=\phi U$ in (\ref{4.23}). Then

\begin{multi} \label{4.25}
g(AU,X)+g(A\phi U,\phi X)=-g(AU,X)
\end{multi}

\noindent for any $X \in {\D}_U$. From (\ref{4.24}) and (\ref{4.25}) it follows

\begin{multi} \label{4.26}
g(AU,X)=g(A\phi U,X)=0
\end{multi}

\noindent for any $X \in {\D}_U$.

The scalar product of (\ref{4.20}) and $U$ yields

\begin{multi} \label{4.27}
g(\phi Y,U)g(A\phi U,X)+kg(\phi Y,U)g(\phi X,U)-g(AX,Y)+g(\phi A\phi Y,X)   \\
=g(Y,U)g(AX,U)-kg(Y,U)g(X,U)    \\
\end{multi}

\noindent for any $X,Y \in \D$. In (\ref{4.27}) we take $Y=U$ and obtain $2g(AU,X)+g(A\phi U,\phi X)=kg(U,X)$, for any $X \in \D$. If $X=\phi U$ we get

\begin{multi} \label{4.28}
g(AU,\phi U)=0
\end{multi}

\noindent and if $X=U$ we have

\begin{multi} \label{4.29}
2g(AU,U)+g(A\phi U,\phi U)=k.
\end{multi}

Taking $Y=\phi U$ in (\ref{4.27}) we obtain $2g(A\phi U,X)+g(\phi AU,X)=kg(\phi U,X)$. If $X=\phi U$ we conclude

\begin{multi} \label{4.30}
2g(A\phi U,\phi U)+g(AU,U)=k.
\end{multi}

From (\ref{4.29}) and (\ref{4.30})

\begin{multi} \label{4.31}
g(AU,U)=g(A\phi U,\phi U)=\frac{k}{3}.
\end{multi}

From (\ref{4.26}), (\ref{4.28}) and (\ref{4.31}) we obtain

\begin{multi} \label{4.32}
AU=\beta\xi +\frac{k}{3} U,      \\
A\phi U=\frac{k}{3} \phi U.    \\
\end{multi}

The scalar product of (\ref{4.20}) and $\xi$ yields

\begin{multi} \label{4.33}
g(\phi AX,A\phi Y)-g(\phi A^2\phi Y,X)-\alpha g(AX,Y)+\alpha g(\phi A\phi Y,X)=0
\end{multi}

\noindent for any $X,Y \in \D$. Taking $X=U$ in (\ref{4.33}) we obtain $g(\phi AU,A\phi Y)+g(A\phi Y,A\phi U)-\alpha g(AU,Y)-\alpha g(A\phi U,\phi Y)=0$, for any $Y \in \D$. From (\ref{4.32}) we get $\left(\frac{2k^2}{9}-\frac{2k\alpha}{3}\right)g(U,Y)=0$, for any $Y \in \D$. Taking $U=Y$,

\begin{multi} \label{4.34}
k=3\alpha.
\end{multi}

If we take $X=\phi U$ in (\ref{4.33}) we have $g(\phi A\phi U,A\phi Y)-g(AU,A\phi Y)-\frac{k\alpha}{3} g(\phi U,Y)+\alpha g(A\phi Y,U)=0$, for any $Y \in \D$. From (\ref{4.32}) and (\ref{4.34}) it follows $\beta^2g(U,\phi Y)=0$ for any $Y \in \D$. If $Y=\phi U$, $\beta =0$, which is impossible and finishes the proof of Theorem 7.

{\bf Remark 2.-} With similar proofs to the ones appearing in this section we could also obtain non-existence results for real hypersurfaces in $\C P^m$, $m \geq 3$, satisfying any of the following conditions:

1) $A_T^{(k)}(\phi X,Y)+\phi A_T^{(k)}(X,Y)=0$, for any $X,Y \in \D$ and any nonnull real number $k$.

2) $A_T^{(k)}(X, \phi Y)+\phi A_T^{(k)}(X,Y)=0$, for any $X,Y \in \D$ and any nonnull real number $k$.

 \hfill $\Box$

{\sc
\begin{trivlist}

\item Juan de Dios P\'{e}rez: jdperez@ugr.es  \\
Departamento de Geometria y Topologia  \\
Universidad de Granada \\
18071 Granada  \\
Spain  \\

\item David P\'{e}rez-L\'{o}pez: davidpl109@correo.ugr.es  \\ 

\end{trivlist}
}
\end{document}